\documentclass[reqno,12pt]{amsart}
\usepackage[a4paper,hmargin={2.5cm,2.5cm},vmargin={2.5cm,2.5cm}]{geometry}
\usepackage{amsthm,amssymb,amsmath,amscd}
\usepackage{amsthm}
\usepackage{pgf, tikz}
\usepackage{graphicx,subcaption}
\usepackage{verbatim}
\usepackage{mathrsfs}
\usepackage{amsrefs}
\usetikzlibrary{calc,decorations.pathmorphing}
\usepackage{enumerate}

\usepackage{enumitem}

\usepackage{xcolor}
\usepackage{hyperref}
\hypersetup{
	colorlinks,
	linkcolor={red!60!black},
	citecolor={green!60!black},
	urlcolor={blue!60!black}
}

\makeatletter
\def\moverlay{\mathpalette\mov@rlay}
\def\mov@rlay#1#2{\leavevmode\vtop{   \baselineskip\z@skip \lineskiplimit-\maxdimen
		\ialign{\hfil$\m@th#1##$\hfil\cr#2\crcr}}}
\newcommand{\charfusion}[3][\mathord]{
	#1{\ifx#1\mathop\vphantom{#2}\fi
		\mathpalette\mov@rlay{#2\cr#3}
	}
	\ifx#1\mathop\expandafter\displaylimits\fi}
\makeatother

\usepackage{todonotes}
\newcommand{\td}[1]{\relax}

\theoremstyle{definition}
\newtheorem{definition}{Definition}%[section]

\theoremstyle{plain}
\newtheorem{theorem}[definition]{{Theorem}}
\newtheorem*{prob}%[definition]
{Open Problems}
{Lemma}
{Theorem}
\newtheorem{lemma}[definition]{{Lemma}}
\newtheorem{cor}[definition]{Corollary}

\newtheorem{fact}[definition]{Fact}
\newtheorem{claim}[definition]{Claim}
{Conjecture}

\theoremstyle{remark}

\numberwithin{equation}{section}
\numberwithin{figure}{section}

\newcommand{\eps}{\varepsilon}

\newcommand{\bv}{\mathbf{v}}

\newcommand{\emm}{\bf}
\newcommand{\bG}{\mathbf{G}}

\begin{document}
		
\title[Size Ramsey numbers for cycles]{On size Ramsey numbers for a pair of cycles}

\author[M.~Bednarska-Bzd\c{e}ga]{Ma{\l}gorzata~Bednarska-Bzd\c{e}ga}
\address{Adam Mickiewicz University, Faculty of Mathematics and Computer Science, Pozna\'n, Poland}
\email{mbed@amu.edu.pl}

\author[T.~{\L}uczak]{Tomasz {\L}uczak}
\address{Adam Mickiewicz University, Faculty of Mathematics and Computer Science, Pozna\'n, Poland}
\email{tomasz@amu.edu.pl}
\thanks{The second author was supported in part by National Science Centre, Poland, grant 2022/47/B/ST1/01517.}
%\thanks{The first author is partially supported by National Science Centre, Poland, grant 2017/27/B/ST1/00873.}

			
\begin{abstract}
We show that there exists an absolute constant $A$ such that the size Ramsey number 
of a pair of cycles $(C_n$, $C_{2d})$, where $4\le 2d\le n$, is bounded from above by $An$.
We also study the restricted size Ramsey number for such a pair. 
\end{abstract}
			
\maketitle

\section{Introduction}

For graphs $G$, $H_1,H_2$  we write $G\to (H_1,H_2)$ if every coloring of the edges of 
$G$ with two colors leads to a monochromatic copy of $H_i$ in the $i$th 
color for some $i=1,2$. 
The {\emm Ramsey number $r(H_1,H_2)$} is defined as the minimum number of vertices 
in a graph $G$ for which $G\to (H_1,H_2)$, the {\emm size Ramsey number} $\hat r(H_1,H_2)$  is the minimum number of edges in~$G$ with such a property, while the 
{\emm restricted size Ramsey number}  {$r^*(H_1,H_2)$} is  
the minimum number of edges in a graph $G$ on $r(H_1,H_2)$ vertices such that $G\to (H_1,H_2)$.

In the paper we study the  Ramsey numbers for 
two cycles, when the shorter cycle is even. 
The Ramsey number for a pair of cycles was computed near half a century  ago by 
 Rosta~\cite{R} and independently  Faudree and Schelp~\cite{FS}. 
 In particular, they proved that if $3d\le n$, then 
\begin{equation}\label{eq:i1}
    r(C_n,C_{2d})=n+d-1\,.
\end{equation}
The size Ramsey number for two cycles is much harder to compute -- as far as we know 
its exact value was found only for few cases, when $n$ and $d$ are small. 
Haxell, Kohayakawa, and {\L}uczak~\cite{HKL} proved that for some absolute constants 
$C$ and $A$ we have $\hat r(C_n,C_{2d})\le An$,
whenever $n\ge 2d\ge C\log n $.  However, 
since their argument relied on a sparse version of Szemeredi's Regularity Lemma, both  constants $A$ and $C$ 
were enormous. This result was improved in several ways. Letzter~\cite{Let} and Ara\'ujo, Pavez-Signe\'e, Sanhueza-Matalama~\cite{A} showed that for every $\eta>0$ there exists a constant $\bar A_\eta$ such that for every $n$ large enough there exists a graph $\bar H_n(\eta)$ 
on $n$ vertices with  fewer than $\bar A_\eta n$ edges with the property that $H_n(\eta)\to (C_n,C_n)$ and moreover $H_n(\eta)$ has fewer than $(1+\eta)r(C_n,C_n)$ vertices. The non-diagonal case was studied  by Javadi, Khoeini,  Omidi, and  Pokrovskiy \cite{JKOP} who showed that
 if $n$ is large enough and $n\ge d\ge \log n + 15$, then 
%\begin{equation}\label{eq:i2}
 $\hat r(C_n, C_{2d})\le 10^{12} n$.
%\end{equation}

Our aim is to prove that a similar result holds also when the even cycle is 
shorter than $\log n$.

\begin{theorem}\label{thm:main}
For every $\eta>0$ there exists a constant $A_\eta$ such that  for every  $n$ and $d$, $n\ge 2d$,  there exists a graph $G_n(\eta)$ with at most $(1+\eta)r(C_n,C_{2d})$ vertices and $A_\eta n$ edges, such that $G_n(\eta)\to (C_n,C_{2d}).$
\end{theorem}

Since the constant $A_\eta$ grows as $\eta\to 0$ it is natural 
to ask what happens when the number of vertices of $G_n(\eta)$
is as small as possible, i.e. what is the value of $r^*(C_n, C_{2d})$. 
Although the notion of the restricted size Ramsey number was defined 
some time ago (see Faudree and Schelp~\cite{FS2}), there are rather few results on this type of Ramsey number in the literature. However, unlike the size Ramsey number, in some cases it can by computed exactly. In particular,  
{\L}uczak, Polcyn, and Rahimi~\cite{LPZ} proved that if $n\ge \ell\ge3 $, $\ell$ is odd, 
and $n$ is large enough, then
\begin{equation*}\label{eq:i3}
 r^*(C_n,C_\ell)=   \lceil (n+1)(2n-1)/2\rceil\,.
\end{equation*}
We remark that the value of $r^*(C_n,C_\ell)$, as well as the value of 
$r(C_n,C_\ell)$, does not depend on the length of the shorter cycle provided it is odd, which makes the analysis somewhat easier. 
This is not  the case when the shorter cycle is even.  The following result gives some estimates for  $r^*(C_n,C_{2d})$.

\begin{theorem}\label{thm:main2}
For every $d\ge 2$ and $n\ge 64d$    
 $$n\max \Big\{\frac{d+1}2, \frac{\log_2 (n/d)}{8\log_2\log_2 (n/d)}\Big\} \le r^*(C_n, C_{2d})\le 20dn \log_2(n/d)\,.
 $$

 Moreover, there exists a constant $d_0$ such that 
 $$r^*(C_n, C_{2d}) \le 20 dn \max\Big\{100+\log_2\frac{n\log_2d}{d^2},1\Big\}$$
for $d\ge d_0$. In particular, for  $d\ge \max\{d_0, 10^{15}\sqrt{n\ln n}\}$ we have
$$r^*(C_n, C_{2d}) \le 20 dn.$$ 
\end{theorem}

We should emphasize that in the paper we try rather to simplify the argument than to get 
best possible constants. Thus, for instance, after some additional work, one can replace constant 20 in 
Theorem~\ref{thm:main2} by, say, 15. Nonetheless, using our argument one cannot improve it to, say, $1/2+o(1)$, so that it will match the lower bound for $r^*(C_n, C_{2d})$ (for more discussion on it see Section~\ref{s:remarks}).

The structure of the paper is the following. In the next section we detour slightly from our main topic and prove the following result on $\hat r(C_n, K_{d,d})$ (a similar result has been independently proved by Sibilska~\cite{mgr}).

\begin{theorem}\label{thm:main3}
There exists an absolute constant $A>0$ such that for every $0<\eta\le 80/A$, $d\ge 8$, and $n\ge 8d/\eta$, 
there exists a graph $G$ with fewer than $(1+\eta)n$ vertices and fewer than 
$20nd/\eta$ edges such that $G\to (C_n, K_{d,d})$. 

In particular,  for every $d\ge 2$ and $n\ge Ad$ we have 
$$nd/2<\hat r(C_n, K_{d,d})\le Adn.$$
\end{theorem}

There are two reasons for including Theorem~\ref{thm:main3} in the paper on Ramsey numbers of two cycles. First, we are going to estimate $\hat r(C_n,C_{2d})$ 
only for  $d>d_1$, for some large constant $d_1$, so in order to deduce from it Theorem~\ref{thm:main} we need to have an upper  for $\hat r (C_n, C_{2d})$ in the form 
$f(d)n$ for some function $f(d)$. Second, the construction crucial for the proof of Theorem~\ref{thm:main} is based on the same idea as the one we use in the proof of Theorem~\ref{thm:main3}, although to verify its validity we need to invoke a sparse version of the Regularity Lemma as it is explained in  
 Section~\ref{s:reg}. Section~\ref{s:lower} gives  lower estimates for the restricted 
 size Ramsey number $r^*(C_n,C_{2d})$, while in Section~\ref{s:upper} we estimate this number from above.   
 We conclude the paper with some remarks and open questions.

\section{Proof of Theorem~\ref{thm:main3}}

Our proof of Theorem~\ref{thm:main3} is, in large part, based on the paper of 
Pokrovskiy and Sudakov~\cite{sud}, who computed the value of $r(C_n,K_{d,d})$.

\begin{theorem}\label{thm:sud} 
If $d\ge 8$ and  $n\ge 2\cdot 10^{49}d$, then $r(C_n,K_{d,d})=n+d-1$.
\end{theorem}

The following result  (see Lemma~3.7 in \cite{sud}) is a crucial ingredient of their argument.

\begin{lemma}%[\cite{sud}]
\label{l:sud}
Let $t$ and $d$ be integers with $t \ge 2\cdot 10^{49}d$ and $d \ge 8$. Let $G$ be a graph 
such that $K_{d,d}\not\subseteq G^c$ and $|S\cup N(S)|\geq t$ for every $S \subseteq  V(G)$ with $|S|\geq  d$. If $x$ and $y$ are two endpoints of a path on at least $8d$ vertices in $G$,
then there is a path on $t$ vertices in $G$ with endpoints $x$ and $y$.
\end{lemma}

Let us rewrite the above lemma in the colored graph setting.

\begin{cor}\label{cor:sud1}
Let $t$ and $d$ be integers with $t \ge 2\cdot 10^{49}d$ and $d \ge 8$. Let $G$ denote a complete graph on at least $r(C_t,K_{d,d})=t+d-1$ vertices and suppose every edge of $G$ is colored either red or blue in such a way that   
there exists no red copy of $K_{d,d}$ in $G$.  
Moreover, let $x$ and $y$ be two endpoints of a blue path on at least $8d$ vertices in $G$. 
Then for every $\ell$ such that $t\le \ell \le |V(G)|-d+1$, there exists a blue path on $\ell$ vertices in $G$ with endpoints $x$ and $y$.
\end{cor}

\begin{proof}[Proof of Theorem~\ref{thm:main3}]
Let $A'\ge 20$ be a constant such that $r(C_k, K_{d,d}))=k+d-1$ for every $d\ge 8$ and $k\ge A' d$. From Theorem~\ref{thm:sud} we know that such a constant exists. 

Let $d\ge 8$, $\eta\le 1/A'$ and $n\ge 8d/\eta$. We put $s= \lfloor \eta n/4d \rfloor$ and consider a `blow-up' $G$ of  the cycle $C_{2s}=v_1v_2\cdots v_{2s}$ in which we  replace each vertex $v_i$ by a  clique $S_i$, such that $12d/\eta\ge |S_i|\ge d/\eta$ and $\sum_{i=1}^{2s}|S_i|=n+2sd$, and moreover each edge of $C_{2s}$ is replaced by a complete bipartite graph between corresponding sets. Observe that such a family of sets exists since
$2s\lceil d/\eta\rceil\le n+2sd\le 2s \lfloor 12d/\eta\rfloor$. 
Then, clearly, $G$ has $n+2sd\le (1+\eta) n$ vertices,  while the number of edges of $G$ is bounded from above by 
$$
\frac12(n+2sd)3\cdot\frac{12d}{\eta}\le \frac{18}{\eta}(1+\eta)dn<  \frac{20}{\eta}dn\,.
$$ 
Now let us suppose that we color edges of $G$ with red and blue 
in such a way that there are no red copy of $K_{d,d}$. Then, since $|S_{2i-1}\cup S_{2i}|\ge 2d/\eta\ge(A'+1)d$, by Theorem~\ref{thm:sud} we infer that for $i=1,2,\dots, s$ the subgraph $H_i=G[S_{2i-1}\cup S_{2i}]$ contains a blue cycle $C^i$ on $|S_{2i-1}|+|S_{2i}|-d$ vertices. Moreover, since subgraphs induced by $S_{2i}\cup S_{2i+1}$ contain no red $K_{d,d}$, from each such subgraph we can choose one blue edge $\{x^i,y ^{i+1}\}$ which connects
$C^i$ to $C^{i+1}$ for $i=1,2,\dots, s$, and a blue edge $\{x^{s}, y^1\}$ joining cycles 
$C^s$ and $C^1$. Note that  $x^i\in S_{2i}$ and $y^i\in S_{2i-1}$ are two different vertices of $C^{i}$ and
since the cycle $C^{i}$ is much longer than $16d$, by 
Corollary~\ref{cor:sud1}, for every $i$ we can connect $x^i$ and $y^i$ by a blue path of length $|S_{2i-1}|+|S_{2i}|-2d$ to build  a blue cycle of length $n$ in $G$. Hence the first part of the assertion follows. 

In order to verify the second part, put $A=80A'$, $\eta=1/A'=80/A$ and notice that for $d\ge 8$ and $n\ge Ad/4>8d/\eta$ the above graph $G$ has less than $20A'dn=Adn/4$ edges. 
Now to complete the proof it is enough to observe that $K_{d,d}\subseteq K_{8,8}$ for every $d\le 8$, so the upper bound in the second part of the theorem holds for  $2\le d\le 8$ as well. To see the lower bound notice that every minimal graph $H$ such that $H\to(C_n,K_{d,d})$ has minimum degree greater than $d$. 
\end{proof}

\section{Proof of Theorem~\ref{thm:main}}\label{s:reg}

The argument we use to show Theorem~\ref{thm:main} is somewhat similar to
the one we apply in the proof of Theorem~\ref{thm:main3}. We split 
the set of $(1+\eta) n$ vertices into small sets, arrange them in the circle, 
and in each consecutive ones we embed the same graph $G$. Now, however, 
$G$ is not the complete graph but a sparse graph $G$ which has 
the property that each coloring of the edges of $G$ with red and blue either leads to a red cycle $C_{2d}$, or results in a large set $S$ such that  
every pair of vertices of $S$  is connected by a blue path of any length chosen from some fairly large interval. Then the argument goes as in the proof of 
Theorem ~\ref{thm:main} -- we generate such sets $S$ in every second pair of sets, connect them by edges, and select blue paths in $S$'s in such a way that the resulting blue cycle has length $n$. Thus, the main challenge here is to find
 a graph $G$ we use in this construction. 

In order to accomplish that we employ 
 %an analysis of random graph decomposition obtained from  
a sparse version of the Regularity Lemma, so let us start with some notions necessary to state it correctly. Let $G=(V,E)$ be a graph and $V_1,V_2\subseteq V$ be two disjoint subsets of its vertices. By $e_G(V_1,V_2)$ we mean the number of edges between $V_1$ and $V_2$, and for  
$p>0$ we define the scaled density $d_{p,G}(V_1,V_2)$ as 
$$d_{p,G}(V_1,V_2)= \frac{e(V_1,V_2)}{p|V_1| |V_2|}\,.$$
Here and below we shall often omit the index $G$ in $e_G(V_1,V_2)$ and $d_{p,G}(V_1,V_2)$
when it does not lead to misunderstandings.
A pair  $(V_1,V_2)$ is called {\emm $(p,\eps)$-regular} if for every $U_i\subseteq V_i$, 
such that $|U_i|\ge \eps|V_i|$, for $i=1,2$, we have 
$$\big|d_p(V_1,V_2)-d_p(U_1,U_2)\big|\le \eps\,.$$
We call such a pair $(V_1,V_2)$  {\emm good} 
if for every $W_i\subseteq V_i$, we have 
$$|N(W_i)|\ge \min\{9 |W_i|, (1-2\eps)|V_{3-i}|\}, \quad \textrm{for}\quad i=1,2,$$
where by $N_G(W) =N(W)$ we always denote the neigborhood of the set $W$ in a graph $G$.
It turns out that, as was proved by  Balogh, Csaba and Samotij~\cite{BCS},  every $(p,\eps)$-regular pair contains a large good $(p,\eps)$-regular pair.
The following lemma is a special case of Lemma~19 from~\cite{BCS}.

\begin{lemma}\label{l:BCS}
 Let  $(V_1,V_2)$ be an $(p,\eps)$-regular  pair for some $0<\eps <0.1$ such that $(1-2\eps)k\le |V_1|,|V_2|\le k$ and $d_p(V_1,V_2)\ge \eps $. 
 Moreover, let $V'_i\subseteq V_i$, $|V'_i|\ge  40 \eps k$, for $i=1,2$. Then there exist sets
 $V''_i\subseteq V'_i$ such that $|V''_i|\ge (1-\eps )|V'_i|$ and $(V''_1, V''_2)$ is a good 
 $(p, 2\eps |V_1|/|V''_1|)$-regular pair.
\end{lemma}

We also define  an {\emm $\alpha$-expanding tree} as a rooted tree $T$ of height $r$ such that the number of vertices $S_i$ at the distance $i$ from the root is $\lceil\alpha |S_{i-1}|\rceil$ for $i= 1,\dots,r-1$. 
The following fact is crucial for our argument.

\begin{lemma}\label{l:path1}
Let  $(V_1,V_2)$ be a  good $(p,\eps)$-regular pair for some $0<\eps \le 10^{-4}$ such that $(1-2\eps)k\le |V_1|,|V_2|\le k$ and $d_p(V_1,V_2)\ge 3\eps$.  Moreover, 
let $h=\lceil\log_8(\eps k)\rceil$,  $h\le \ell\le (2-10\sqrt {\eps})k$, and $x\in V_1\cup V_2$.
Then there exists a set $Y\subseteq V_1\cup V_2$ of at least $\eps k$ vertices such that for each $y\in Y$ there is a path of length $\ell$ joining  $x$ and $y$.
\end{lemma}

\begin{proof}  We shall show, using an induction on $\ell$, a slightly stronger statement, namely, that for every $\ell$,  $h\le \ell\le (2-10\sqrt {\eps})k$, and   $x\in V_1\cup V_2$, 
there exist sets $Z$ and $Y\subseteq Z$ such that $|Z|\le \ell +4\eps k$, $Y\ge \eps k$, 
and for each $y\in Y$ there exists a path $P_{x,y}$ of length $\ell$ joining  $x$ and $y$ such 
that all vertices of $P_{x,y}$ are contained in $Z$.

Note first that since 
the pair $(V_1,V_2)$ is good,
for every $2\le \alpha\le 8$ and every  $x\in V_1\cup V_2$, the vertex $x$ is a root of an $\alpha$-expanding tree $T(x,\alpha)$ such that the number of leaves at the highest level of $T(x,\alpha)$ is $\lceil 2\eps k\rceil$. In particular, the assertion holds for all $\ell$
such that $h\le \ell\le 3h$.

Now let us suppose that the assertion holds for some $x$ and $\ell_0$, and by  $Z_0$ and $Y_0$ 
we denote the sets which certify that it is true. 
We shall show that it holds also for $x$ and  $\ell_1=\ell_0+h$. 
In order to see it consider the pair $(V'_1,V'_2)$, 
where $V_i\setminus Z_0$ for $i=1,2$.
From Lemma~\ref{l:BCS} it follows that  there exist sets $V''_i\subseteq V'_i$, $i=1,2$, such that $|V''_1|, |V''_2|\ge \sqrt{\eps}k\ge \eps k$ and the pair $(V''_1,V''_2)$ is a good $(p,2\sqrt{\eps})$-regular pair. Since $(V_1,V_2)$ was $(p,\eps)$-regular, 
there exists at least one edge $e=\{y,z\}$ between the set $Y_0$ and $V''_1\cup V''_2$; let 
$y=e\cap Y_0$ and 
$z=e\cap (V''_1\cup V''_2)$. However, since $(V''_1,V''_2)$ is good, there exists a
4-expanding tree $T(z,4)$ rooted at $z$ which has at least 
$$4\sqrt{\eps}\min\{|V''_1|,|V''_2|\}\ge \eps k$$
leaves. Thus, we may take for $Z$ the vertices of a path 
of length $\ell_0$ joining  $x$ and $y$ and all vertices of $T(z,4)$, and for $Y$ 
 the set of vertices in distance $\ell_1$ from $z$ in $T(z,4)$. 
\end{proof}

Let us note the following consequence of the above lemma.

\begin{lemma}\label{l:path2}
Let  $(V_1,V_2)$ be a  good $(p,\eps)$-regular pair for some $0<\eps <10^{-4}$ such that $(1-\eps)k\le |V_1|,|V_2|\le k$ 
and $d_p(V_1,V_2)\ge 3\eps$. 
Moreover, 
let $h=\lceil\log_4(\eps k)\rceil$,  $\ell$ be an odd natural number such that  $4h\le \ell\le (2-30\sqrt{\eps})k$, and $x_i\in V_i$ 
for $i=1,2$.
Then there exists a path of length $\ell$ joining  $x_1$ and $x_2$.

In particular, $(V_1,V_2)$ contains a cycle of length $\ell +1$.
\end{lemma}

\begin{proof}  
Using the expanding property we build two vertex disjoint 4-expanding trees $T_1(x_1,4)$ and $T_2(x_2,4)$ 
of height $h$ rooted at $x_1$ and $x_2$ respectively. Let 
$V'_1=V_1\setminus (V(T_1)\cup V(T_2))$ and $V'_2=V_2\setminus (V(T_1)\cup V(T_2))$.
From Lemma~\ref{l:BCS} it follows that there exist sets $V''_i\subseteq V'_i$, $i=1,2$, such that $|V''_i|\ge (1-8\sqrt {\eps})k$
and $(V''_1,V''_2)$ is a good $(p, 2\sqrt{\eps})$-regular pair. 
Let $x\in V''_1\cup V''_2$ be one of 
neighbors of at least $\eps k$
leaves of $T_1(x_1,4)$ -- since the pair $(V_1,V_2)$ 
is $(p,\eps)$-regular such a neighbor always exists. Using Lemma~\ref{l:path1}
we generate a set $Y$ of at least $\eps k $ vertices each of which is connected to $x$ by 
a path of length $\ell-2h-2$. Then, using again  $(p,\eps)$-regularity of
$(V_1,V_2)$, we find an edge connecting one of the vertices of $Y$ with one of the leaves of $T_2(4,\eps)$, closing a path of length $\ell$ between $x_1$ and $x_2$. 

To see  the last part of the lemma it is enough to take any edge and join its ends by a path of length $\ell$. 
\end{proof}

Let us recall that a sparse version of the Regularity Lemma,  discovered independently by Kohayakawa and R\"odl
in the early nineties, states that for every $\eps>0$ 
the vertex set of a graph $G=(V,E)$ whose edges are, in a way, `uniformly distributed', 
can be partitioned into a few sets of equal size so that almost all pairs of sets are 
$(p, \eps)$-regular, where $p$ is the density $|E|/\binom{|V|}2$ of $G$. 
The condition that edges are `uniformly distributed' holds for  random graphs and all its 
dense subgraphs, so the Regularity Lemma can be used.  Since this applications of the Regularity Lemma is now routine (see, for instance,  \cites{A, BCS, HKL, Let}) we only state here its consequence, when it is applied to a random graph $G(N,c/N)$. 

\begin{lemma}\label{l:reg}
For every $\eps>0$  there exist constants $N_0$, $T$, and $c$, such that the following holds.
For every $N\ge N_0$ there exists a graph $G=(V,E)$, $|V|=N$, $|E|\le cN/2$, with the property that  for every coloring of edges of $G$ with red and blue there exists a partition of $V$ into sets $V_1,\dots, V_t$, such that
\begin{enumerate}
    \item [(i)] $1/\eps\le t\le T$;
    \item [(ii)]  $\big||V_i|-|V_j|\big |\le 1$ for $1\le i<j\le t$;
    \item[(iii)] all but at most $\eps t^2$ pairs  $(V_i,V_j)$ are $(c/N,\eps)$-regular in both the red graph $R$ and the blue graph $B$;
    \item[(iv)] for every $1\le i<j\le t$ and $p=c/N$ either  $d_{p,B}(V_i,V_j)\ge 1/3$ or $d_{p,R}(V_i,V_j)\ge 1/3$.
\end{enumerate}
\end{lemma}

We denote the graph whose existence is guaranteed by the above lemma by  $\hat G_N(c,\eps)$.
For any coloring of  $\hat G_N(c,\eps)$ and a partition for which conditions (i)--(iv) 
hold, by $\bG_t(\eps)$ we denote the {\emm reduced graph} of the partition defined as 
the graph with 
vertices $\bv_1,\dots, \bv_t$, where two vertices $\bv_i$ and $\bv_j$ are connected by a red [blue]
edge if the pair $(V_i,V_j)$ is $(p,\eps)$-regular in  the red graph $R$ [the blue graph $B$] with $p=c/N$ and its 
scaled density in $R$ [$B$] is larger than $1/3$. Note that $\bG_t(\eps)$ is the complete graph on $t$ vertices from which we have removed at most $\eps t^2$ edges and remaining edges are colored with red and blue. 

The following result is, in some way, analogous to Lemma~\ref{l:sud}. 

\begin{lemma}\label{l:path3}
For every positive $\eta<1$ there exist constants $n_1$, $a_1$, $A_1$ and $c_1$ such that 
for every $n\ge n_1$ there exists a graph $G_n(\eta)$ with  $n(1+\eta)$ vertices and 
fewer than $c_1n/2$ edges for which the following property holds.

For every $d$ and $\ell$ such that  $A_1\log n \le d\le a_1 n$, $A_1\log n\le \ell\le n(1+\eta/3)$,  and every coloring of edges of $G_n(\eta)$ with red and blue which do not lead to red $C_{2d}$:
\begin{enumerate}
    \item [(i)] each two subsets of vertices of $G_n(\eta)$ of size $n/3$ each are connected by at least one blue edge;
    \item[(ii)] 
there exists a vertex set $S$ such that $|S|\ge (1+\eta/2)n$ and any two vertices $x,y\in S$
are connected by a blue path of length $\ell$. 
\end{enumerate}

In particular, for each such coloring, the graph
 $G_n(\eta)$ contains a blue cycle of length $\ell+1$. 
\end{lemma}

\begin{proof}
Put $\eps=\eta^2/10^4$, $N=(1+\eta)n$, and let $\hat G_N(c,\eps)$ be a graph whose existence is 
assured by Lemma~\ref{l:reg}. Color its edges with red and blue so that it contains 
no red copy of $C_{2d}$. First we argue that the reduced graph $\bG_t(\eps)$ contains 
no red edges. Indeed, such edge means that a red graph contains an $(p,\eps)$-regular pair 
$(V_i,V_j)$ with scaled density at least 1/3 and $\lfloor n/t\rfloor \le 
|V_i|, |V_j|\le \lceil N/t\rceil $ and so, by Lemma~\ref{l:path2}, it contains a (red) cycle of every length between $\log_4 (N/t)$ and $n/(2t)$. Hence, $\bG_t(\eps)$ is a graph on $t$
vertices with at least $(1-2\eps)\binom t2$ blue edges. 
Note that  such a graph can have at most 
$3\eps t$ vertices of degree smaller than $t/3$. It proves (i) and, 
due to  Dirac Theorem, it implies that the reduced graph contains  an odd blue cycle on 
at least $2r+1\ge (1-4\eps)t$ vertices. Without loss of generality let us assume that this is
a cycle $\bv_1\bv_2\cdots \bv_{2r+1}\bv_1$. 
Now let $W^j_i\subseteq V_i$, $|W^j_i|=\lfloor \sqrt\eps N/t\rfloor$, for $i=1,2,\dots, 2r$ and $j=1,2,3,4$.  Using Lemma~\ref{l:BCS} 
we infer that  for every $i=1,2,\dots 2r$, there exists $\bar V_i\subseteq V_i\setminus (W^1_i\cup W^2_i \cup W^3_i\cup W^4_i)$ such that
$|\bar V_j|\ge (1-5\sqrt\eps)N/t$ and for $s=1,2,\dots, r$
the pair $(\bar V_{2s-1},\bar V_{2s})$ is a good  $(p,\eps)$-regular pair, with $p$-density at least $1/3$ for $p=c/N$. Moreover, 
since for every $s=1,2,\dots, r$, and $j=1,2,3,4$, the pair $(W_{2s-1}^j, W_{2s-1}^j)$ is $(p,\sqrt\eps)$-regular, there exist sets 
$\bar W^j_{2s-1}\subseteq W^j_{2s-1}$ with $|\bar W^j_{2s-1}|\ge (1-2\sqrt{\eps})|W^j_{2s-1}|$, and 
$\bar W^j_{2s}\subseteq W^j_{2s}$ with $|\bar W^j_{2s}|\ge (1-2\sqrt{\eps})|W^j_{2s}|$, such that pairs $(\bar W_{2s-1}^j, \bar W_{2s}^j)$ are
good $(p,2\sqrt\eps)$-regular pairs. 
Note also that since  the pair $(V_{2s-1},V_{2s})$ is $(p,\eps)$-regular with scaled density at least 1/3,
there exist sets  $\hat V_{2s-1}\subseteq \bar V_{2s-1}$,  and  $\hat V_{2s}\subseteq \bar V_{2s}$, such that 
$|\hat V_{2s-1}|\ge (1-6\eps)|\bar V_{2s-1}|$, $|\hat V_{2s}|\ge (1-6\eps)|\bar V_{2s}|$,
 each vertex of $\hat V_{2s-1}$ has at least 
one neighbor in $\hat V_{2s}$ and each of the sets $\bar W^j_{2s}$, and 
 each vertex of $\hat V_{2s}$ has at least 
one neighbor in $\hat V_{2s-1}$ and in each of the sets $\bar W^j_{2s-1}$, for $s=1,2,\dots, r$ and   
for $j=1,2,3,4$. 

Let us set 
$$S=\bigcup_{i=1}^{2r}\hat V_i\,.$$
Note first that 
$$|S|\ge 2r (1-6\eps)(1-5\sqrt\eps) N/t\ge \big((1-4\eps)t -1\big)(1-6\sqrt\eps)N/t\ge (1-7\sqrt\eps)N\ge (1+\eta/2)n\,.$$
Now let  $x,y\in S$. We will argue that in $\hat G_N(c,\eps)$ there are even and odd blue paths joining $x$ and $y$, both  of length at most 
$13r\log_4(N/t)$, which contain precisely one edge from each pair $(\bar V_{2s-1}, \bar V_{2s})$ for $s=1,2,\dots, r$.  
In order to construct such paths we use the fact that 
$\bv_1\bv_2\cdots \bv_{2r+1}\bv_1$
is a blue cycle in  the reduced graph $\bG_t(\eps)$. Let us suppose that $x\in \hat V_1$ and $y\in \hat V_{i_0}$ for some $i_0=1,2,\dots, 2r$. 
First we connect  $x$ to a vertex $v_2$ from  $\hat V_2$. For the next vertex of the path we choose a neighbor $w_1\in \bar W^1_1$ of $v_2$. 
Then we use the fact that the pair $(\bar W^1_1, \bar W_2^1)$  is good and build a tree $T_1$ rooted at $w_1$ with $\eps |V_2|$ leaves in $\bar W^1_2$. 
Since the pair $(V_2,V_3)$ is $(p,\eps) $-regular, at least one of these leaves has a neighbor $w_3$ in $\bar W^1_3$ -- we select it as the next vertex of the path.  
Then we use edges of the pair $(\bar W^1_3, \bar W^1_4)$ to build a tree $T_2$ rooted at $w_3$ which has 
at least $\eps |V_3|$ leaves in, say, $\bar W^1_4$. At least one of these leaves has 
 a  neighbor $v_3$ in  $\hat V_3$; we select it as the next vertex of our path and 
 connect it with one of its neighbors $v_4\in \hat V_4$. Further, we choose 
a vertex  $w'_3$ in $\bar W^2_3$ adjacent to $v_4$ and build a tree rooted on $w'_3$ with a lot of leaves, using edges between $(\bar W^2_3,\bar W^2_4)$. 
In this way  we go around the cycle using `buffer sets' $\bar W^j_i$ and $V_{2r+1}$,
picking one edge from each pair $(\hat V_{2s-1}, \hat V_{2s})$ on the way, but never use any buffer more than ones. Moreover, we omit the pair $(\hat V_{2s-1}, \hat V_{2s})$ which contains $y$ until the very end.
To close the path  we start with $y$, 
choose its neighbor $z$ in $\hat V_{i_0+1}$ or $\hat V_{i_0-1}$,  choose a neighbor $z_{i_0}$ of $z$ in 
buffer set $\bar W^j_{i_0}$  which has not being used so far, for some $j=1,2,3,4$, and build a tree with a lot 
of leaves rooted at $z_{i_0}$, so that we can connect it to rest of the path. Note also that, 
if necessary, we can adjust the parity of the path running over the cycle, using buffer sets, one extra time. 
Thus, one can create a path joining $x$ and $y$ 
%\TL{mysle ze teraz jest juz dobrze}
of required parity and length at most $4(2r+1)\log_4 N\le 13r \log_4 n$, which contains one edge from every pair $(\bar V_{2s-1}, \bar V_{2s})$ and only one vertex from each of the sets
$\bar V_i$, $i=1,2,\dots,2r$. Now we can apply  Lemma~\ref{l:path2} to
adjust the length of the path by replacing each blue edge contained in a good pair 
by a path of any odd length up to 
$$(2-30\sqrt{\eps})(1-4\sqrt\eps)N/t\ge (2-40\sqrt{\eps})N/t\,.$$
Thus $x$ and $y$ can be joined by a path of any length $\ell$, provided 
$\ell\ge 13r\log_4(N/t)$ and 
$$\ell\le r(2-40\sqrt{\eps})\frac{N}{t}\,. $$
Since 
$$  r(2-40\sqrt{\eps})\frac{N}{t}\ge  \frac{(1-4\eps)t}{2}(2-40\sqrt{\eps})\frac{N}{t}\ge (1-21\sqrt{\eps})(1+\eta)n\ge (1+\eta/3)n\,,   $$
this completes the proof of (ii).
\end{proof}

\begin{proof} [Proof of Theorem~\ref{thm:main}]
 If $d$ is smaller than $n/C$ and larger than $C\log n$ for some large constant 
 $C$, the assertion follows from Lemma~\ref{l:path3}. If $d$ is of the order $n$,
 then $r(C_n,C_{2d})$ is  larger than $n(1+\eta)$, but  the argument similar to that 
 we used to prove Lemma~\ref{l:path3} still works -- one should find in the reduced graph 
 either large enough blue or large enough red cycle and apply Lemma~\ref{l:path2} to adjust it length. 
 Since it is a rather standard procedure
 (see Letzter~\cite{Let2} for recent developments in this area) and the proof 
 follows the very same line as in the diagonal case studied in \cites{A, Let},
 we omit the details  and concentrate on the most interesting case, when 
 $d=O(\log n)$. 
 
We show Theorem~\ref{thm:main} for $d$ which is large enough, i.e. larger than  $d_1=\sqrt{n_1}$, where $n_1$ is a constant chosen in such a way the the assertion of  Lemma~\ref{l:path3} holds. Our argument will be somewhat similar to the construction we used in 
the proof of Theorem~\ref{thm:main3}. For $0<\eta<1/4$, $n$, and $d$ such that  $d_1\le d\le n^{1/3}$, let us consider a graph $G_{d^2}(\eta)$
on  $N=\lfloor (1+\eta)d^2\rfloor $ vertices and average degree smaller than some constant $c_1$, 
 so that for $G_{d^2}(\eta)$ the assertion of Lemma~\ref{l:path3} holds. 
Now let us build a new graph 
$H=(V,E)$ whose vertex set consists of disjoint sets $W_1,\dots, W_{2r}$, where
$|W_{2i}|=\lfloor N/2\rfloor$, $|W_{2i-1}|=\lceil N/2\rceil$,
for $i=1,2,\dots, r$ and $r=\lfloor n/d^2 \rfloor$. In every pair $W_i\cup  W_{i+1}$, $i=1,2,\dots, 2r$,
we embed a copy of graph $G_{d^2}(\eta)$ and we also embed such a copy in a pair 
$W_1\cup W_{2r}$. Note that $H$ has 
$$r(1+\eta)d^2\le (1+\eta)n$$
vertices and the average degree of $H$ is bounded from above by $2c_1$.
Now let us suppose that we color edges of  $H$ with red and blue in such a way that there are no red copies of $C_{2d}$. Then in view of Lemma~\ref{l:path3} each copy of $G_{d^2}(\eta)$ 
embedded in $W_{2i-1}\cup W_{2i}$ for $i=1,\dots, r$, contains a set $S_i$ such that
$|S_i|\ge (1+\eta/2)d^2$ and  each two vertices of $S_i$ are connected by path of any length $\ell$ between,
say, $d^2/2$ and $(1+\eta/3)d^2$. Moreover, $S_i$ and $S_{i+1}$ 
are connected by at least one blue edge for $i=1,2,\dots, r-1$, and so are the sets $S_1$ and $S_r$. 
Consequently, the graph $H$ contains blue cycles of any length $\ell$
such that
$$r d^2/2 \le \ell \le r (1+\eta/3)d^2 \,.$$
Since 
$$rd^2/2\le n(1+\eta)/2\le n$$
and 
$$r (1+\eta/3)d^2\ge (1+\eta/4)n\ge n,$$
$H$ contains a blue copy of $C_n$. This completes the proof for $d\ge d_1$. 
However, from Theorem~\ref{thm:main3} we know that for $d\le d_1$
$$\hat r(C_n, C_{2d})\le \hat r(C_n, K_{d,d})\le Ad n\,,$$
so the assertion of Theorem~\ref{thm:main} holds for small $d$ as well. 
\end{proof}

\section{Lower bounds for $r^*(C_n, C_{2d})$}\label{s:lower}

Since for  $n\ge 3d$ we have  $r(C_n,C_{2d})=n+d-1$, we 
 start  with the following observation.

\begin{fact}\label{f:1}
If $G\to (C_n,C_{2d})$ and $G$ has $n+d-1$ vertices, then $\delta(G)\ge d+1$. 

In particular, if $n\ge 3d$, then 
$r^*(C_n, C_{2d})\ge \lceil n(d+1)/2\rceil$.
\end{fact}

\begin{proof}
Let us suppose that some vertex  $v$ of $G$ has neighbors 
$w_1,w_2,\dots, w_d$ (if the degree 
of $v$ is smaller than $d$, just add to $G$ some edges). Now color all edges incident to one of the vertices $w_1,w_2,\dots, w_{d-1}$, 
as well as the edge $\{v,w_d\}$, with red, and all 
the remaining edges with blue. In this coloring no  $C_n$ is colored blue and no  $C_{2d}$ is colored red. 
\end{proof}

 When $d$ is small, this bound can be substantially improved. Here and below $P_n$ denote the path on $n$ vertices.

\begin{theorem}\label{thm:lower}
Let $n\ge 64$ and $1\le b\le n/64$. Suppose that 
$H=(V,E)$ is a graph on $n+b-1$ vertices such that $H\to (P_4,P_n)$.
Then 
$$|E|\ge\frac{n\log_2(n/b)}{8\log_2\log_2(n/b)}\,.$$ 
\end{theorem}

\begin{proof}
Let $s=4|E|/|V|$. 
We define recursively two graph sequences $H_0\supseteq H_1\ldots\supseteq H_t$ and $G_0\supseteq G_1\supseteq \ldots \supseteq  G_t$, a non-decreasing sequence 
$X_0 \subseteq X_1\subseteq \ldots\subseteq X_t$ 
of  subsets of $V$   and 
one more sequence $(S_1,S_2,\ldots, S_t)$ of subsets of $V$ in the following way. Put
$X_0=\emptyset$, $H_0=H$, and let $G_0$ be the graph induced in $H$ on the set of all vertices of degree at most $s$ in $H$. Let $S_j$ be the largest set of vertices of $G_j$ such that their pairwise distances in $H_j$ are greater than 2. 
As long as $N_{H_{j-1}}(S_j)$ is not empty, we define 
the set $X_j=X_{j-1}\cup N_{H_{j-1}}(S_j)$ and the graphs $H_j=H_{j-1}\setminus  N_{H_{j-1}}(S_j)$, $G_j=G_{j-1}\setminus  N_{H_{j-1}}(S_j)$.  
By $t$ we mean the largest index  $j\ge 1$ such that $N_{H_{j-1}}(S_{j})$ is not empty.

\begin{claim}\label{neighb1}
For every $0\le j\le t$ the following holds.
\begin{enumerate}
\item[(i)]
$|S_{j+1}|\le |X_j|+b$, 
\item[(ii)]
$|N_{H_{j}}(S_{j+1})|\le (|X_j|+b)s$,
\item[(iii)] 
$V(G_j)\subseteq S_{j+1}\cup N_{H_j}(S_{j+1})\cup N_{G_j}(N_{H_j}(S_{j+1}))$.
\end{enumerate}
\end{claim}

\begin{proof}
Suppose that $|S_{j+1}|> |X_{j}|+b$,  for some $0\le j\le  t$.
Observe that if we color the edges of $H$ such that all edges incident to $S_{j+1}$ in $H_j$ are red and all remaining edges of $H$ are blue, then we have no red copy of $P_4$, as well as no blue path longer than $|V(H_j)|-|S_{j+1}|+2|X_{j}|+1\le |V(H)|-b=n-1$. It contradicts the assumption that $H\to(P_4,P_n)$ and hence (i) follows. 
 
The second part of the claim is the consequence of the first one and the fact that every vertex in $S_{j+1}\subseteq V(G_j)$ has degree at most $s$.

To see (iii) note that if there exists a vertex $u\in V(G_j)\setminus S_{j+1}$ such that 
$u\not\in N_{H_j}(S_{j+1})\cup N_{G_j}(N_{H_j}(S_{j+1}))$, then $u$ is in distance at least 3 from every vertex of $S_{j+1}$, 
which contradicts the maximality of $S_{j+1}$.
 \end{proof}

Based on the above claim and the definition of $X_i$, one can easily verify inductively that for every $0\le j\le t$ we have 
\begin{equation}\label{neighb2}
|N_{H_{j}}(S_{j+1})|\le (s+1)^j bs\text{ and }
|X_j|\le \big((s+1)^j-1\big)b.
\end{equation}
Furthermore, since $V(G_j)\subseteq S_{j+1}\cup N_{H_{j}}(S_{j+1})\cup N_{G_{j}}(N_{H_{j}}(S_{j+1}))$, every vertex of $V(G_j)\setminus N_{H_j}(S_{j+1})$ either is isolated in $H_j$ or has got at least one neighbor in $N_{H_{j}}(S_{j+1})$. Therefore
$\text{deg}_{H_{j+1}}(u)\le \text{deg}_{H_j}(u)-1$ for every vertex $u\in V(G_{j+1})$ with  positive degree in $H_j$. Thus
 $t\le s$, since the degree in $H$ of every vertex of $u\in V(G_0)$ is at most $s$.  
 
Let us recall that we assumed that $N_{H_{t}}(S_{t+1})$ is empty. It means that $S_{t+1}=V(G_t)$. The construction of the sets $X_j$ implies that $|V(G_0)|\le |V(G_t)|+|X_t|$ so  $|V(G_0)|\le |S_{t+1}|+|X_t|$ and by Claim \ref{neighb1} we have $|V(G_0)|\le 2|X_t|+b$.
Hence by (\ref{neighb2}) we obtain  
$$|V(G_0)|\le 2\big((s+1)^t-1\big)b+b\le 2\big((s+1)^s-1\big)b+b<2(s+1)^s b.$$

In view of the definition of $s$ and $G_0$, we have $|V(G_0)|\ge |V|/2$, so the above inequality implies that 
$(s+1)^s>|V|/(4b)$ and hence 
$s\ge \log_2(|V|/(4b))/\log_2\log_2(|V|/(4b))$ for $b\le |V|/64$. 
Therefore  
$$|E|=\frac14 |V|s\ge \frac{|V| \log_2(|V|/(4b))}{4\log_2\log_2(|V|/(4b))}\ge \frac{n\log_2(n/b)}{8\log_2\log_2(n/b)}\,.$$ 
\end{proof}

The above result  immediately implies the following lower 
bound for $r^*(C_n, C_{2d})$, which for  $d\ll \log n/\log\log n$ is better than 
the simple bound given in Fact~\ref{f:1}, but, unfortunately, it decreases 
as a function of $d$.

\begin{cor}
 If $n\ge 64d$, then 
$$r^*(C_n, C_{2d})\ge \frac{n \log_2(n/d)}{8\log_2\log_2(n/d)}.$$ 
\end{cor}

We also remark that it follows from Theorem~\ref{thm:lower}  that 
$r^*(C_n, P_4)=\Omega(n\log n/\log\log n)$, while it is known that 
$r^*(C_n,P_3)=O(n)$ (see for instance Ben-Shimon, Krivelevich, Sudakov~\cite{BKS}).

%\MB{Czy nie dobrze byloby wspomniec juz w Introduction, ze jako wynik uboczny mamy cos  o $r*(P_4,P_n)$? Czasem czytelnicy (tacy jak ja) szukajac wynikow na jakis temat, czytaja tylko intro, a potem okazuje sie, ze autorzy ukryli w glebi pracy to, czego czytelnik szuka.}

\section{Upper bounds for $r^*(C_n, C_{2d})$}\label{s:upper}

In order to estimate $r^*(C_n, C_{2d})$ from above we need to find a sparse 
graph $G$ on $n+d-1$ vertices such that $G\to (C_n, C_{2d})$.
We start with a technical lemma which says the  following.  
Let $F$ be a (sparse) graph such that for every coloring of its edges 
with red and blue which does not lead to a red copy of $C_{2d}$, we find 
at most $s$ blue paths which covers all vertices of $H$. Then, if we add 
 $O(s+d)$ vertices to $F$ and connect them with every vertex of $F$, then for the
resulting graph $G$ we have $G\to (C_n,C_{2d})$, where $|V(G)|=n-d+1$, precisely 
as we want.  In order to state this result   we 
introduce an additional notion. 
We call a partially colored graph $\hat{G}$
on $n$ vertices 
an {\emm $(n,s,t)$-system}, if $\hat G$ consists of:
\begin{itemize}
\item 
the `central clique' $K$ on $t$ vertices, 
\item
 $s\ge 0$ vertex disjoint `satellite paths' $L_1,\ldots, L_s$, 
whose edges are colored with blue. All vertices of each of the  paths $L_i$, $i=1,2,\dots, s$,  are connected by 
(uncolored) edges to the central clique.
\end{itemize}

Note that some of the satellite paths may be  trivial, i.e. they may be  
just vertices.

\begin{lemma}\label{mainlemma}
Let $G$ be an $(n,s,t)$-system  such that  $t\ge 10d+4s$ for some $d\ge 2$. 
Then each coloring of the uncolored edges of $G$ with two colors, red 
and blue, leads to either a red copy of $C_{2d}$, or a blue copy of $C_{n-d+1}$.
\end{lemma}

\begin{proof}
Let $K$ be the central clique of $G$ and $L_1,L_2,\ldots,L_s$ be the family of the satellite paths.
We denote the fully colored $G$ by  $\bar G$ and assume that it contains 
no red copy of $C_{2d}$.
By  $X$ we denote the  set of $d-1$  vertices of $\bar G$ which have the smallest blue neighborhood 
in the central clique and  write  $G'$ for the subgraph induced in $\bar G$ by $V(G)\setminus X$. 

\begin{claim}\label{bigblue}
Every vertex of $V(G')$ has at least $4d+2s$ blue neighbors in the central clique. 
\end{claim}

\begin{proof} 
Suppose that a vertex  $u\in V(G')$ has fewer  than $4d+2s$ blue neighbors in the central clique. 
Let $A=X\cup\{u\}$ and $B\subseteq V(K\setminus X)$.
Then $|A|=d$, $|B|\le |V(K)|-|X|>t-d\ge 9d+4s$ 
and so every two vertices in $A$ have more than $|B|-8d-4s\ge d$ common red neighbors in $B$.
Hence in a greedy way we can find a red $C_{2d}$ in $\bar G[A\cup B]$.
\end{proof}

%Put $J=\{i:\, V(L_i)\setminus X\neq \emptyset\}$. 
For every $i=1,2,\dots,s$,  we define a path $L'_i$ as 
the longest path contained in $L_i$ such that its ends do not belong to $X$ (if all vertices of $L_i$ are in $X$  we say that the path is  empty so we can ignore it). Note 
that $V(L'_i)\supseteq V(L_i)\setminus X$ but $V(L'_i)$ can also contain some vertices
from $X$. Let $d'$ stand for  the number of such  vertices, i.e. 
$$d'=\Big|X\cap \bigcup_{i=1}^s V(L'_i)\Big|\,.$$ 
We `compensate' for these vertices by selecting in $V(K)\setminus X$ any $d'$ vertices
and denote this set as $X'$. 
Note that
$$|V(K)\cap (X\cup X')|=|V(K)\cap X|+d' \le d-1\,.$$
We will show that the blue subgraph spanned in $\bar G$ by $(V(K)\setminus X')\cup \bigcup_{i=1}^s V(L'_i)$ is hamiltonian. 

Let us first study graph  $H$ which is  the blue subgraph induced in $\bar G$ by vertices $V(K)\setminus (X\cup X')$.

\begin{claim}
$H$  is $2(d+s)$-vertex-connected. 
\end{claim}

\begin{proof} 
By the previous claim, each vertex of $V(H)$ has at least $4d+2s$ blue neighbors in $K$ and hence it has at least $4d+2s-|V(K)\cap (X\cup X')|>3d+2s$ blue neighbors in $H$. Moreover,  since $K$ contains no red copy of $C_{2d}$, any pair of disjoint  subsets $U_1, U_2\subseteq V(H)$ of $d$ vertices each is connected by at least one blue edge. Thus, every vertex cut of $H$ has at least $2d+2s$ vertices.  
\end{proof}

We also remark  that, in view of Claim \ref{bigblue} and the definition of $H$, every vertex  of $V(G')$ 
has at least $4d+2s-|V(K)\cap (X\cup X')|> 3d+2s$ blue neighbors in $H$.

%\TL{Malgosiu, w akapicie pod spodem zmienilem "incident to all vertices" na "incident to at least $3d+2s$ vertices". Czy slusznie?}

For every $i=1,2,\dots,s$ such that $L'_i$ is not  empty 
we select a set $E_i$ of blue edges in the following way. 
By the definition of the path $L'_i$, its ends $x,y$ (possibly equal) have more than $3d+2s$ blue neighbors in $V(H)$.
%incident to more at least $3d+2s$ vertices of the central clique. 
We choose blue edges $xa_i, yb_i\in E(\bar G)$ such that $a_i,b_i\in V(H)$, $a_i\neq b_i$
and we define the set $E_i=E(L'_i)\cup\{xa_i,yb_i\}$. Note that $E_i$ is the edge set of a blue path. 
Moreover, since we choose blue neighbors of  at most  $2s$ vertices, and each such  vertex  
has more than $3d+2s$ blue neighbors in $V(H)$, we can select $a_i $'s and $b_i$'s in such a way 
that they are all distinct. 

Now we add at most $2s$ auxiliary edges  $\{a_i,b_i\}$  to $H$ and color them azure (if $a_ib_i\in E(H)$, then we repaint it). 
In this way we obtain the blue-azure graph $H'$ consisting of all non-repainted blue edges of $H$ and a matching $B$ which consists of 
at most $s$ azure edges. Clearly, for the connectivity $\kappa(H')$ of $H'$ we get
$$\kappa(H')\ge \kappa (H)\ge 2(d+s)\,.$$ 
Furthermore, from the assumption that $\bar G$ has no red $C_{2d}$, we infer that 
the independence number  $\alpha(H')$ of $H'$ is smaller than $2d$.
Hence,  $\kappa(H')>\alpha(H')+|B|$. Now we use the following result of 
 H\"aggkvist and Thomassen~\cite{thom} which is a generalization of the well known
 Chv\'atal-Erd\H{o}s criterion for hamiltonicity. 

\begin{lemma}%[\cite{thom}]
\label{ech}
If $H'$ is a graph with $\kappa(H')\ge \alpha(H')+m$, then for every set of vertex disjoint paths in $H'$ with $m$ 
edges in total, $H'$ has a Hamilton cycle containing all these paths.  
\end{lemma}

Thus, the graph $H'$ contains a Hamilton cycle $C'$ containing all azure edges. 
It is not hard to see that  the edge set $(E(C')\setminus B)\cup \bigcup_{i=1}^s E_i$ forms a blue cycle $C$ in $\bar G$ such that
$V(G)\setminus V(C)$ consists of all elements of $X\setminus\bigcup_{i=1}^s V(L'_i)$ and all elements of $V(K)\cap (X\cup X')$. Thus, in view of the definition of $X'$, exactly $d-1$ vertices of $G$ do not belong to $C$. Thereby we obtain a blue cycle on $n-d+1$ vertices. 
\end{proof}
 
We use the above result to estimate $r^*(C_n,C_{2d})$ for $d=\Omega(\sqrt{n\ln n})$.

\begin{lemma}\label{l:dlarge}
If $10^{13}\sqrt{n\ln n}\le 2d\le n/10$ and $d$ is large enough, then 
\begin{equation}
r^*(C_n,C_{2d})\le     20dn\,.
\end{equation}
\end{lemma}

\begin{proof}
The following observation is crucial for our argument.

\begin{claim}
For every large enough $d$ such that $10^{13}\sqrt{n\ln n}\le 2d\le n$ 
there exists a graph $G(n,d)$ with $N=n-17d-1$ vertices and fewer than $dn$ edges
such that each subset of vertices of $G(n,d)$ on at least $4d$ vertices contains 
a copy of $C_{2d}$.
\end{claim}
\begin{proof}
Let us consider a random graph $G(N,p)$, where $p=d/3n$. Its edges are binomially distributed with parameters $\binom N2<n^2/2$ and $p$, so 
with probability at least 1/3 it has fewer than $dn$ edges. 
Now let $\eps=10^{-4}$ and let random variable $X$ counts 
pairs of disjoint sets  of vertices $S,T$  of 
$G(N,p)$ such that $|S|, |T|\ge 2\eps d $ and the number of edges
 deviates from its expected value $p|S||T| $ by more than $\eps p|S||T|$.
Then, using Chernoff bounds (see, for instance, Corollary~2.3 in \cite{RG})
\begin{align*}
 \textsc{E} X&\le \sum_{i=\lceil 2\eps d\rceil }^N \sum_{j=\lceil 2\eps d\rceil }^N \binom {N}{i} \binom {N}{j} \exp\Big(-\frac{\eps^2}3 ijp \Big)
 \le N^2 \binom {N}{\lceil 2\eps d\rceil } ^2 
 \exp\Big(-\frac{\eps^2 (2\eps d)^2 d}{9N} \Big)\\
 &\le 
 N^2\Big(\frac{e^2N^2}{4\eps^2 d^2}\exp\Big(-\frac{2\eps^3d^2}{3N} \Big)\Big)^{2\eps d}, 
\end{align*} 
and,  since $d\ge 10^{13}\sqrt{n\ln n}$, $N\le n$ and $2\cdot 10^{13}\eps^{3}/9\ge 2$,
for $d$ large enough we get
$$\Pr(X> 0)\le \textrm{E} X\le n^2 \big(e^2\eps^{-2}/d \big)^{2\eps d}\le 1/3\,.$$
Thus, with probability at least $2/3$ each pair of subsets $(U,W)$ 
with $|U|, |W|\ge 2d$, is a $(p,\eps)$-regular pair.  However, 
due to Lemmata~\ref{l:BCS} and~\ref{l:path1}, such a pair contains 
cycles of each length from, say, $d/4$, to $2d(2-2\eps -10\sqrt{\eps})>2d$.
Hence, with positive probability  a graph $G(n,d)$ with required property exists.
\end{proof}
To build a graph $H(n,d)$ with $n+d-1$ vertices and fewer than $20dn$ edges
such that $H(n,d)\to (C_n,C_{2d})$, take 
a graph $G(n,d)$ whose existence is assured in Claim above, add to it a clique on $18d$ new vertices, 
and connect every vertex of the clique to all vertices of $G(n,d)$. 
Now let us color the edges of a graph $H(n,d)$ obtained in this way 
with red and blue
so that there are no red copies of $C_{2d}$. Consider the largest 
family $\mathcal {P}$ of vertex disjoint blue paths contained in $G(n,d)$, and 
let $S$ denote the set obtained by taking  one end of each such path. 
From the maximality of the family $\mathcal {P}$ we infer 
that $S$ contains no blue edges, and so $\mathcal {P}$ contains fewer than $4d$ 
paths. Now we can use Lemma~\ref{mainlemma} with $s=4d$ to deduce that
in this coloring of $H(n,d)$ there exists a blue copy of~$C_n$.
\end{proof} 

In order to present a construction which gives a general upper bound  
for $r^*(C_n,C_{2d})$, we introduce some more definition. By $T_N$ 
we denote a {\emm binary tree} on $N$ vertices which is obtained from the perfect
rooted binary tree of height $h=\lceil \log_2 (N- 1)\rceil $ by removing some 
leaves on the highest level. Let $\hat T_N$ be the {\emm closure} of $T_N$
which is obtained from $\hat T_N$ by joining each vertex of $T_N$ with each of its
descendants. By leaves of $\hat T_n$ we mean the vertices which had degree 1 in 
$T_N$ (and have  degree  $h -1$ or $h-2$ in $\hat T_N$).   
Note  also that if $N=2^k-1$, then vertices of each level of $\hat T_N$ send 
at most $N$ edges to its descendants, so the number of edges in $\hat T_N$ 
can be crudely estimated from above by $N\log_2 N$. 

Now for $n\ge 14d \ge 28$  let $U(n,d)$ be a 'blow-up' of the closure of a binary tree, 
which is constructed in the following way. 
Take  $N=\lfloor (n+d-1)/(14d) \rfloor$ and replace each vertex of $\hat T_N$
by a clique of $14d$ elements, except, perhaps,  one leaf which we replace by the clique
of 
$$(n+d-1)-14d(N-1)\le 28d$$
 elements, so that the resulting graph has precisely $n+d-1$ vertices. 
 Moreover, we replace each edge by the complete bipartite graph between the corresponding sets.
Thus, $U(n,d)$ has at most 
$$(14d)^2N\log_2 N+N\binom{14d}{2}+(14d) (\log_2 n+14d) \le 20dn\log_2(n/d)$$
edges, provided $n\ge 14d$.

The following observation is a consequence of Lemma~\ref{mainlemma}.

 \begin{lemma}\label{0treearrow}
If  $n\ge 14d\ge 28$, then $U(n,d)\to (C_n, C_{2d})$.  

In particular, for $n\ge 3d\ge 6$ we have $r^*(C_n,C_{2d})\le 20nd \log_2(n/d)$. 
\end{lemma}

\begin{proof}
We show that $U(n,d)\to (C_n, C_{2d})$ using induction on $n$. 
If $n< 28d$, then $U(n,d)$ is a clique of size larger than $r(C_n,C_{2d})=n+d-1$, so the assertion holds.
Let us assume that $n\ge 28d$ and let $K$ denote the set of $14d$ vertices which replaced the root of $\hat T_N$.
Then, if we remove $K$ from $U(n,d)$, the resulting graph $\bar U(n,d)$ either  can be identified with $U(n-14d, d)$
or consists of two components which are graphs $U(n_1,d)$ and $U(n_2,d)$, where
 $n_1+n_2+d-1=n-14d$. Now suppose that we color edges of $U(n,d)$ by two colors, red and blue,
 in such a way that there are no red copies of $C_{2d}$.
Then, by the induction hypothesis, in the obtained graph $\bar U(n,d)$ there exists a blue path which contains all but at most $d-1$
vertices of the graph, or there exist two blue paths which contain all but at most $2(d-1)$ vertices combined. Since we treat isolated vertices  as 
trivial paths, it means that the vertex set of $\bar U(n,d)$ can be covered by at most $2d$ blue paths. Since the clique $K$ has $14d$ vertices, we infer from Lemma~\ref{mainlemma}
 that one can use
the edges of $K$ and the edges between $K$ and $\bar U(n,d)$ to create a blue cycle on precisely $n$ vertices.

To see the second part of Lemma~\ref{0treearrow} it is enough to observe that if $2d\le n\le 14d$, then the complete graph 
on $n+d-1$ vertices has fewer than $20dn$ edges.
\end{proof}

Notice that one can easily construct graphs $G$ on $n+d-1$ vertices such that $G\to (C_n, C_{2d})$ and have density slightly smaller than $U(n,d)$. For instance, 
in  the closure of an appropriately chosen binary tree we may replace all leaves not by 
 cliques of size $14d$, but by independent sets of size $d$.
%\TL{teraz chyba jest lepiej}
However, this improvement only modifies a constant next to $dn$ which we, quite crudely, estimated by 20. However, when $d$ is close to $\sqrt n$ one can use Lemma~\ref{l:dlarge}
and  can get a substantially better estimates by replacing each leaf by a copy of a graph $G$ 
on $k+d-1=10^{-26}d^2/\log_2 d$ vertices and at most $20dk$ edges which is  such that $G\to (C_k, C_{2d})$.  Since $k\gg d$, the height of a tree we need is  $O(\log_2 ((n+d-1)/k))$ so it results in $O(dn\log (n/k))$ as the upper bound for $r^*(C_n, C_{2d})$ 
(Theorem~\ref{thm:main2} gives bounds with specific constants).  

\section{Final remarks}\label{s:remarks}

Clearly, in the paper we do not resolve all problems concerning 
the size Ramsey and the restricted size Ramsey numbers for pair of cycles.  The first question which naturally comes to mind is whether $\hat r(C_n, C_{\ell})=O(n) $
also in the case when $\ell\le n$ and $\ell $ is odd. We are convinced 
that this is the case and believe that we can prove it, however, since the argument 
is much more complicated than in the case of even $\ell$, we decided not to include it 
here, especially that this work is dedicated to study Ramsey numbers of pairs of cycles in which  the shorter cycle is even. 

Theorem~\ref{thm:main2} raises much more questions. Let us state 
at least some of them. Here by $o_d(1)$  we denote a function which tends to 0
 as $d\to\infty$.

\begin{prob}
Find smallest possible constants $0\le a_1\le a_2\le a_3$ such that for every $0<\eta< 1-a_3$ and $3d\le n$
large enough the following holds.
\begin{enumerate}
    \item [(i)] For $2d\ge n^{a_3+\eta}$ we have 
    \begin{equation}\label{eq:last1} 
    r^*(C_n, C_{2d})= \lceil (n+d-1)(d+1)/2\rceil\,.
    \end{equation}
     \item [(ii)] For $2d\ge n^{a_2+\eta}$ we have 
    \begin{equation}\label{eq:last2} 
    r^*(C_n, C_{2d})= (1/2+o_d(1))d(n+d)\,.
    \end{equation}
     \item [(iii)] For $2d\ge n^{a_1+\eta}$ we have
    \begin{equation}\label{eq:last3} 
    r^*(C_n, C_{2d})\le 20 dn\,.
    \end{equation}
\end{enumerate}
\end{prob}

Theorem~\ref{thm:main2} implies that $a_1\le 1/2$. On the other hand, it is easy to see 
that $a_3\ge 1/2$. Indeed, each graph $G$ on $n+d-1$ vertices with $\Delta(G)\le \sqrt n/2$
contains a vertex $v_1$ and vertices $v_2,\dots,v_d\notin N(N(v_1))$. Then we can color 
all edges incident to one of vertices $v_1,v_2,\dots,v_d$ red and the rest of edges 
blue, creating neither red copy of $C_{2d}$ nor blue copy of $C_n$. Hence 
$G\not\to (C_n, C_{2d})$. 

Besides of the above information  we can only speculate on the values $a_1,a_2,a_3$.
Our guess, which is not based on any solid evidence, is that $a_3=2/3$, $a_2=1/2$, 
and the value of $a_1$ is either $0$ or $1/2$. 

Finally, another interesting problem is to study the asymptotic behavior of $r^*(C_n,C_4)$ 
and $r^*(P_n,P_4)$ which we determined only up to a factor of $\log\log n$.

\end{document}